\documentclass[11pt,a4]{article}

\usepackage[T2A]{fontenc}
\usepackage[utf8]{inputenc}
\usepackage[english,russian]{babel}
\usepackage{xcolor}
\usepackage{indentfirst}
\usepackage{graphicx}
\usepackage{epstopdf}
\usepackage{eqnarray,amsmath}
\usepackage{amsfonts, amsmath, amsthm}
\usepackage{amssymb}
\usepackage{url}
\usepackage{caption}
\usepackage{subcaption}
\def\ai#1{{\color{black}#1}} 
\def\ag#1{{\color{black}#1}} 
\usepackage{algorithm}
\usepackage{algpseudocode}
\usepackage{graphicx}
\usepackage{kbordermatrix}

\def\algorithmicif{\textbf{если}}
\def\algorithmicthen{\textbf{то}}

\def\algorithmicdпo{}

\newcommand{\argmin}{\mathop{\arg\,\min}}
\newcommand{\argmax}{\mathop{\arg\,\max}}
\DeclareMathOperator*{\argminset}{Arg\,\min}

\newtheorem{Theo}{\bf ~~~~Теорема}

\newtheorem{Lem}{\bf ~~~~Лемма}
\newtheorem{Coro}{\bf ~~~~Утверждение}

\begin{document}

\begin{center}
\textbf{Поиск равновесия по Вальрасу и централизованная распределённая оптимизация с точки зрения современных численных методов выпуклой оптимизации на примере задачи распределения ресурсов}
\end{center}

\begin{center}
\textit{Евгения Алексеевна Воронцова
\footnote{vorontsovaea@gmail.com, Дальневосточный федеральный университет, 690091, Россия, г. Владивосток, ул. Суханова, 8\\
Evgeniya Vorontsova, Far Eastern Federal University, 8 Sukhanova St.,
Vladivostok, Russia 690091
}
Александр Владимирович Гасников
\footnote{gasnikov.av@mipt.ru, Московский физико-технический институт (государственный университет), 141700, Россия, Московская область, г. Долгопрудный, Институтский переулок, д. 9; Институт проблем передачи информации им. А.А. Харкевича Российской академии наук, 127051, Россия, г. Москва, Большой Каретный переулок, д.19 стр. 1; Адыгейский государственный университет, 385000, Россия, Майкоп, ул. Первомайская, 208 \\
Alexander Gasnikov, Moscow Institute of Physics and Technology, 9 Institutskiy per., Dolgoprudny, Moscow Region,
Russia 141700; Institute for Information Transmission Problems RAS, Bolshoy Karetny per. 19, build.1,
Moscow, Russia 127051; Adyghe State University, Str. Day, 208, Maikop, Republic of Adygea, 385000}
Анастасия Сергеевна Иванова
\footnote{anastasiya.s.ivanova@phystech.edu,  Московский физико-технический институт (государственный университет), 141700, Россия, Московская область, г. Долгопрудный, Институтский переулок, д. 9\\
Anastasiya Ivanova, Moscow Institute of Physics and Technology, 9 Institutskiy per., Dolgoprudny, Moscow Region,
Russia 141700
}
Евгений Алексеевич Нурминский
\footnote{nurmi@dvo.ru, Дальневосточный федеральный университет, 690091, Россия, г. Владивосток, ул. Суханова, 8\\
Evgeny Nurminsky, Far Eastern Federal University, 8 Sukhanova St.,
Vladivostok, Russia 690091
}}
\end{center}

\begin{abstract}
В данной работе на примере численного решения классической задачи распределения ресурсов демонстрируются: 1) Вальрасов механизм нащупывания равновесия; 2) Децентрализующая роль цен; 3)  Слейтеровская конструкция по ограничению цен (двойственных множителей); 4) Новый механизм поиска равновесных цен, в котором цены устанавливаются не 
Центром (государством), а узлами (предприятиями). В отличие от экономической литературы, в которой, в основном, ограничиваются установлением факта сходимости исследуемых процедур, в работе приводится точный анализ скорости сходимости описываемых процедур поиска равновесия с учётом их прямо-двойственной природы. По сути, в  работе предпринята попытка содержательно (экономически) проинтерпретировать следующие численные процедуры одновременного решения прямых и двойственных задач выпуклой оптимизации: метод дихотомии и метод проекции субградиента.
\end{abstract}

\textbf{Keywords} Walrasian equilibrium, decentralized pricing, primal-dual method, subgradient method, Slater condition

\textbf{Ключевые слова} 
Вальрасов механизм, децентрализация цен, прямо-двойственный метод, субградиентный метод, условие Слейтера

\textbf{Mathematics Subject Classification}  90B99	

\textbf{УДК} 519.86

\section*{Введение}

В работе рассматривается задача \ag{(в максимально упрощенной формулировке)} распределения ресурсов в централизованно и децентрализованно (распределёно) управляемой экономике \cite[глава 10]{Oben}. Предлагаются новые способы интерпретации решения данной задачи.

Статья организована следующим образом. 

В п.~1 рассматривается классическая задача распределения ресурсов (глава 10 \cite{Oben}). Центр (государство) является собственником предприятий,
которые производят один и тот же товар. Центру известны функции денежных затрат всех предприятий на производство данного товара. Задача Центра
 заключается в том, чтобы так спланировать производственные задания (планы) для подведомственных предприятий, что при заданном объёме производства суммарные затраты предприятий будут минимальными.  
 Такая задача может быть решена с помощью стандартных приёмов математического программирования \cite{Boyd}, \cite{Nedic}. 

В п.~2 рассматривается модификация данной задачи, в которой Центр уже не является собственником предприятий (у каждого предприятия есть свой собственник). Как следствие, функции затрат предприятий Центру, вообще говоря, неизвестны. Центр устанавливает цену на товар, производимый всеми предприятиями
(модель по-прежнему однопродуктовая). Исходя из этой цены, предприятия (оптимально для себя) производят товар. В зависимости от объёма произведённого товара Центр выставляет новую цену товара на следующий этап (год). Такой механизм <<нащупывания>> равновесной цены в экономической литературе называют вальрасовым механизмом \cite[глава 10]{Oben}. В данной работе в основу вальрасова механизма был положен метод дихотомии. Насколько нам известно, ранее на вальрасов механизм в основном смотрели как на локальную процедуру корректировки цены \ag{(исключением является работа \cite{Friedman})}. В данной статье предпринята попытка уйти от этого ограничения. Как следствие, удалось получить линейную скорость сходимости, которая для рассматриваемого класса задач будет оптимальной. Процедура дихотомии предполагает ограниченность возможного диапазона цены товара. Начальная область локализации цены определяется из условий Слейтера. Таким образом, в п.~2 построен численный алгоритм решения задачи распределения ресурсов. Алгоритм содержательно проинтерпретирован. Получена не улучшаемая с точностью до числового множителя для рассматриваемого класса задач оценка скорости сходимости. В этом отличие результатов п.~2 от большинства экономических работ, предлагающих только
доказательство факта сходимости предлагаемых процедур нащупывания равновесия. 

В п.~3 приведена интерпретация задачи нащупывания цен из п.~2
как задачи централизованной распределённой оптимизации \cite{Parallel and distributed optimization}, \cite{Nedic}.

В п.~4 приведён новый механизм нащупывания цен (развивающий идеи статьи \cite{Primal-dual subgradient methods}), в котором цена устанавливается не Центром (государством), а узлами (предприятиями). Каждое предприятие производит товар и устанавливает
цену для продажи товара Центру. Центр отбирает наиболее выгодные предложения и \ag{пытается закупить} товар в требуемом объёме. Предприятия корректируют объём производимого товара и цены, исходя из того, закупал ли Центр у них товар, и если закупал, то в каком объёме \ag{(еще точнее в зависимости от спроса на товар со стороны Центра)}. Описан алгоритм получения \ag{(нащупывания)} равновесных цен с помощью метода проекции субградиента. Доказана сходимость метода и получена оценка скорости сходимости.

В п.~5 приведены численные эксперименты по моделям из пп.~2,~4.

\ag{Подчеркнем, что в данной работе рассматривается карикатурная модель взаимодействия Центра с предприятиями. В частности, в п. 2  предполагается, что 1) предприятия производят нескладируемый товар (нельзя делать запасы и продавать их в будущем); 2) в случае, если предприятия произведут меньше, чем требуется Центру (имеет место дефицит товара: спрос больше предложения), Центр просто делает выводы и на следующий этап выставляет большую цену закупки; 3) независимо от объема производства Центр покупает у предприятий все, что они произведут. Если первое предположение не выглядит сильно ограничительным, то второе и третье нуждаются в пояснении. Второе предположение, на наш взгляд, вполне может иметь место, и описанная в п. 2 процедура дихотомического подбора равновесной цены позволяет быстро исправить эту ситуацию. Третье предположение является, пожалуй, наименее очевидным. Но и здесь вполне можно представить ситуацию, в которой Центр (например, государство) стремится таким образом поддерживать производителей, обещая им, что все что они смогут произвести будет выкуплено по установленным заранее ценам. Здесь также, как и по второму предложению в случае перепроизводства Центр просто учитывает это, уменьшая закупочную цену на следующий этап. В случае интерпретации модели из п. 4, пожалуй, самым слабым местом является, что Центр закупает товар только у тех предприятий, которые производят товар по самой низкой цене. И даже если суммарный объем производства этих предприятий меньше, чем спрос, то Центр не докупает товар у предприятий, выставивших следующую цену. Такое предположение может соответствовать аукционной модели. В будущем мы постараемся обобщить предложенную здесь модель таким образом, чтобы Центр мог закупать товар и у предприятий, выставивших не самую низкую цену.}

\section{Задача распределения ресурсов при известных функциях затрат}\label{sect_1}
Рассмотрим следующую задачу распределения ресурсов.
Пусть Центр владеет $n$ предприятиями, производящими один и тот же продукт \ag{(товар)}. У каждого предприятия есть своя функция затрат $f_k(x)$, $k = 1, \, \ldots, \, n$, которая зависит от объёма выпуска этого продукта в год $x$ (тонн). Пусть Центру необходимо, чтобы суммарно все предприятия производили не менее чем~$C$ тонн продукта. Кроме того, Центр, как владелец предприятий, хочет минимизировать суммарные затраты на производство, что приводит к следующей задаче:


\ag{
\begin{equation}
\label{2_eq_main}
f(x) = \sum \limits_{k=1}^n f_k(x_k) \rightarrow \min \limits_{\substack{\sum \limits_{k=1}^n x_k \, \geqslant  \, C; \;  \\ 
x_k \, \geqslant  \, 0,\; k =1, \, \ldots, \, n.}}
\end{equation}}
Будем предполагать, что функции затрат $f_k(x)$, $k =1, \, \ldots, \, n$ возрастающие и $\mu$-сильно выпуклые \cite{Oben} (сильная выпуклость обеспечивает существование и единственность решения задач~\eqref{2_eq_main} и \eqref{2_eq_x_k(p)}).

Задача~\eqref{2_eq_main} является частным случаем задач
монотропического программирования (см.~\cite{Bertsekas_mono}, \cite{Rok_mono}).
Применим для решения задачи~\eqref{2_eq_main} метод множителей Лагранжа \cite{Boyd}. Для этого запишем двойственную задачу:

\begin{equation*}
\begin{split}
\min \limits_{\substack{\sum \limits_{k=1}^n x_k \, \geqslant  \, C; \;  \\ 
x_k \, \geqslant  \, 0,\; k =1, \, \ldots, \, n}}\sum \limits_{k=1}^n \;f_k(x_k)& = \min \limits_{x_k \geqslant  0,\; k =1, \, \ldots, \, n}\left\{\sum \limits_{k=1}^n f_k(x_k)+\max \limits_{p \geqslant  0}\left\{\left(C-\sum \limits_{k=1}^n x_k\right)p\right\}\right\}\\
& = \max \limits_{p \geqslant  0}\Bigg\{\sum \limits_{k=1}^n \min \limits_{x_k \geqslant  0} \Big\{f_k(x_k)-px_k\Big\}+pC\Bigg\}\\
& =-\min \limits_{p \geqslant  0}\Bigg\{\sum \limits_{k=1}^n \max \limits_{x_k \geqslant  0} \Big\{px_k-f_k(x_k)\Big\}-pC\Bigg\}\\
& =- \min \limits_{p \geqslant  0}\left\{\sum \limits_{k=1}^n \left\{px_k(p)-f_k(x_k(p))\right\}-pC\right\},
\end{split} 
\end{equation*}
где 
\begin{equation}
\label{2_eq_x_k(p)}
x_k(p)=\argmax \limits_{x_k \geqslant  0} \Big \{ p x_k-f_k(x_k) \Big \}, \quad  k =1, \, \ldots, \, n.
\end{equation}
Тогда двойственная задача (с точностью до знака) имеет следующий вид

\begin{equation}
\label{2_eq_dual}
\varphi(p)=\sum \limits_{k=1}^n \Big\{px_k(p)-f_k(x_k(p))\Big\}-pC \rightarrow \min \limits_{p \geqslant  0}.
\end{equation}

\begin{Coro} \cite{Nest}
\label{2_cor_sopr}
Пусть для любого $k = 1,...,n$  функция $f_k(x)=f_k(x_k)$ является $\mu$-сильно выпуклой.
Тогда функция (\ref{2_eq_dual}), где $x_k(p),\; k =1, \, \ldots, \, n$ определяются условием~\eqref{2_eq_x_k(p)}, является $n/\mu$-гладкой. То есть, производная функции $\varphi(p)$ удовлетворяет условию Липшица с константой $L_{\varphi}= n/\mu$
\begin{equation}
\label{2_smooth}
\left|\varphi'(p_2)-\varphi'(p_1)\right| \leqslant   L_{\varphi}\left|p_2-p_1\right|.
\end{equation}
\end{Coro}

\begin{Coro}~\textbf{(Теорема Демьянова--Данскина--Рубинова)}
\cite{Danskin, Dem_book_74}
\label{2_cor_dem}
Пусть для любого $p\geqslant 0$ выполняется: $\varphi(p)=\max \limits_{x \in X} F(x, \, p)$, $F(x, \, p)$~--- выпуклая и гладкая по $p$ функция и максимум достигается в единственной точке $x(p)$. Тогда $\varphi'(p) = F_{p}'(x(p), \,p)$.
\end{Coro}
В нашем случае для задачи~\eqref{2_eq_dual} из утверждения 2 получаем
\begin{equation}
\label{2_eq_grad_phi}
\varphi'(p)=\sum \limits_{k=1}^n x_k(p)-C.
\end{equation}
Пусть~$p^*$~--- решение двойственной задачи~\eqref{2_eq_dual}. Тогда, согласно необходимому условию экстремума, для $p^*$ должно выполняться условие дополняющей нежесткости
\begin{equation*}
\label{2_eq_extr_1}
p^*\left(\sum \limits_{k=1}^n x_k(p^*)-C\right)=0,\; p^*\geqslant  0,
\end{equation*}
что, учитывая соотношение~\eqref{2_eq_grad_phi}, можно представить в следующем виде 
\begin{equation}
\label{2_eq_extr_2}
p^*\varphi'(p^*)=0,\; p^*\geqslant  0.
\end{equation}
Однако, в силу экономической интерпретации задачи, по которой~$p$ выступает в роли цены, стоит рассматривать только положительные значения для~$p$. Получается, что условие~\eqref{2_eq_extr_2} равносильно следующему условию 

\begin{equation}
\label{2_eq_extr}
\varphi'(p^*)=\sum \limits_{k=1}^n x_k(p^*)-C=0.
\end{equation}

Объём производства для каждого предприятия Центр определяет с учётом желания минимизировать затраты. То есть, устанавливая цену на товар равной~$p$, Центр для каждого предприятия решает следующую задачу нахождения объёма выпуска

\begin{equation}
\label{2_eq_1}
x_k(p)=\argmin \limits_{x_k \geqslant  0} \Big \{ f_k(x_k)-p x_k\Big \}, \quad  k =1, \, \ldots, \, n.
\end{equation}
Сформулируем следующую теорему (доказательство стандартно и базируется на теореме Каруша--Куна--Таккера \cite{Boyd}):

\begin{Theo}
\label{theo_1}
Пусть решением задачи~\eqref{2_eq_main} является вектор~$\mathtt{x^{*}}=(x^{*}_{1}, \, \ldots, \, x^{*}_{n})^{\top}$. Тогда существует~$p^*$ такое, что~$x^*_{k}= x_k(p^*) = \argmin\limits_{x_{k} \, \geqslant  \, 0} \big\{f_{k}(x_{k})-p^*x_{k}\big\}, \; \;  k =1, \, \ldots, \, n$,
где цена~$p^*$ является решением двойственной задачи \eqref{2_eq_dual}, т.е. определяется из условия~\eqref{2_eq_extr}.
\end{Theo}

\section{Задача распределения ресурсов \\
при отсутствии информации о функциях затрат}
\label{sect_no_f}

В предыдущем пункте рассмотрен случай известных Центру функций затрат~--- Центр являлся собственником всех предприятий. В данном пункте рассматривается чаще встречающаяся на рынке ситуация~--- у каждого предприятия имеется свой собственник. Однако Центру по-прежнему необходимо, чтобы предприятия суммарно за год произвели не менее чем~$C$ тонн товара. 

В текущей задаче функции затрат каждого предприятия~$f_k(x_k)$, $k =1, \, \ldots, \, n$ неизвестны Центру,  и, соответственно, 
невозможно воспользоваться
методом решения из п.~\ref{sect_1}. Тем не менее, если 
Центр может устанавливать цену на товар, то посредством этого он может регулировать суммарный объём производства предприятий. Происходит это следующим образом. Центром устанавливается цена~$p$ на товар и сообщается предприятиям. После этого каждое предприятие анализирует, сколько тонн товара оно произведет при такой цене, и сообщает эту информацию Центру. Затем Центр сравнивает суммарный объём производства за год~$\sum \limits_{k=1}^n x_k(p)$ с необходимым количеством товара~$C$, и, в зависимости от результата, корректирует цену.
При этом каждое предприятие определяет объём производства, исходя из желания максимизировать прибыль при заданной цене,
т.е. каждое предприятие решает следующую задачу нахождения объёма выпуска 
\begin{equation*}
\label{3_eq_max_xkp}
x_{k}(p)=\argmax \limits_{x_{k}\geqslant  0} \big\{\underbrace{px_{k}}_\text{выручка}-\underbrace{f_{k}(x_{k})}_\text{затраты}\big\}, \; \;  k =1, \, \ldots, \, n.
\end{equation*}
Общая задача совпадает с~\eqref{2_eq_main}, однако в силу отсутствия у Центра информации о функциях затрат предприятий, решить ее аналитически уже невозможно.
Поэтому для решения задачи предлагается 
использовать один из простейших численных методов~---
метод дихотомии (деления отрезка пополам). Далее будет более подробно описан алгоритм поведения Центра и предприятий
в этом случае.
Также будет получен положительный (в виду теоремы~\ref{theo_1}) ответ на вопрос: если~$p^*$~--- равновесная цена~\eqref{2_eq_extr}, установленная Центром, а~$\mathtt{x}^{*}=(x^{*}_{1}, \, \ldots, \, x^{*}_{n})^{\top}$~---  решение исходной задачи~\eqref{2_eq_main},
будет ли выполняться равенство~$x_k(p^*)=x_k^*, \; \;  k =1, \, \ldots, \, n$?

\subsection{Стратегия Центра по нащупыванию равновесной цены.\\ Поиск равновесия по Вальрасу}
\label{sect_no_f_1}
Пусть Центр каким-то образом локализовал цену в отрезке~$\Delta = [p_{min}, \, p_{max}]$. Тогда на следующем шаге он устанавливает цену~$p=\dfrac{p_{min}+p_{max}}{2}$ и сообщает её предприятиям. После этого каждое предприятие решает свою задачу максимизации прибыли и сообщает Центру, какое количество товара оно произведёт при данной цене.
Получив <<обратную связь>> от предприятий, Центр определяет, произвели ли суммарно предприятия необходимое ему количество товара (т.е., определяет знак~$\varphi'(p)=\sum \limits_{k=1}^n x_k(p)-C$) и, в зависимости от этого, корректирует цену следующим образом:

\begin{itemize}
 \item Если суммарный объём производства предприятий больше, чем необходимо\\ ($\varphi'(p) > 0$), то Центр сдвигает правую границу области локализации цены так, что теперь~$p_{max} = \dfrac{p_{min}+p_{max}}{2}$. В этом случае цена слишком большая и можно суммарно произвести необходимое количество и при более низкой цене.

 \item Если суммарный объём производства предприятий меньше, чем необходимо \\($\varphi'(p) < 0$), то Центр сдвигает левую границу области локализации цены по формуле~$p_{min} = \dfrac{p_{min}+p_{max}}{2}$. В этом случае для того, чтобы получить необходимое количество товара, цену необходимо повысить.

 \item Если~$\varphi'(p) = 0$, то предприятия произвели необходимое количество товара и решение найдено.
\end{itemize}
С каждой итерацией (в нашей интерпретации, с каждым годом) отрезок локализации цены уменьшается вдвое. Тогда очевидно, что после~$N$ таких итераций Центр сократит длину начального промежутка локализации цены в~$2^N$ раз. В качестве критерия останова алгоритма выберем точность решения~$\varepsilon$, т. е. алгоритм завершает работу при
выполнении условия
$\left|C-\sum\limits_{k=1}^{n}x_k(p^{N})\right|\leqslant   \varepsilon$, что, учитывая условие дополняющей нежёсткости~\eqref{2_eq_extr}, эквивалентно  $\left|\sum\limits_{k=1}^{n}x_k(p^*)-\sum\limits_{k=1}^{n}x_k(p^{N})\right|\leqslant   \varepsilon$, здесь $N$~--- номер последней итерации.

Формализованное описание стратегии Центра и предприятий по установлению равновесной цены приведено в алгоритме~\ref{alg_dih}.

\floatname{algorithm}{Алгоритм}
	\begin{algorithm}
		\caption{Стратегия Центра по нащупыванию равновесной цены}\label{alg_dih}
		\begin{algorithmic}[1]
			
			\Require $f_k(x), \; \; k=1, \, \ldots, \, n$~--- сильно выпуклые функции затрат для каждого предприятия; $[p_{min}; \, p_{max}]$~--- начальный отрезок локализации цены.
			\Ensure цена $p^N$, для которой выполняется $\left|C-\sum \limits_{k=1}^n x_k(p^N)\right|\leqslant   \varepsilon$.
			\State $N := 0;$
			\While{ $\left|C-\sum           \limits_{k=1}^n x_k(p)\right| > \varepsilon$}
			\State $p := \dfrac{p_{min}+p_{max}}{2};$
			\State     $x_{k}(p)=\argmax \limits_{x_{k}\geqslant  0} \big\{px_{k}-f_{k}(x_{k})\big\}, \; \; k=1, \, \ldots, \, n;$
			
			\State $\varphi'(p)=\sum \limits_{k=1}^{n}x_k(p)-C;$
			\State \algorithmicif{} $\varphi'(p)>0$, \algorithmicthen{} 
			$p_{max} := \dfrac{p_{min}+p_{max}}{2}$;
			\State \algorithmicif{} $\varphi'(p)<0$, \algorithmicthen{} 
			$p_{min} := \dfrac{p_{min}+p_{max}}{2}$;
			\State $N := N + 1$;
			\EndWhile
			\State $p^N=p;$
			\State\Return $p^N$
		\end{algorithmic}
	\end{algorithm}

Несколько слов о том, как определить в начальный момент область локализации цены~$\Delta$. Очевидно, что в качестве левой границы можно взять~$p_{min}=0$. Это следует из экономической интерпретации задачи, по которой цена не может быть отрицательной, и из того что~$p$~--- двойственный множитель к ограничению типа неравенства. Оценка для~$p_{max}$ будет получена далее.

\subsection{Условие Слейтера}

Необходимо получить оценку значения $p_{max}$. Для этого воспользуемся условием Слейтера и проведём рассуждения, аналогичные \cite{Bertsekas}. Положим~$\bar x_1 = \ldots = \bar x_n = 2C/n$. 
Заметим, что в точке $\mathtt{\bar x} = (\bar x_1, \ldots , \bar x_n)^{\top}$ выполняется условие Слейтера, так как
$$C-\sum \limits_{k=1}^n \bar x_k=C-2C=-C < 0.$$

Поскольку функции затрат $f_k(x_k)$, $k = 1, \, \ldots, \, n$, являются возрастающими
(исходя из экономических соображений), то выполняется
неравенство

\begin{equation*}
\begin{split}
\sum \limits_{k=1}^n f_k(0) & = \min \limits_{x_k \geqslant  0,\; k=1, \, \ldots, \, n}\sum \limits_{k=1}^n f_k(x_k)\\
&=\min \limits_{x_k \geqslant  0,\; k=1, \, \ldots, \, n}\left\{\sum \limits_{k=1}^n f_k(x_k)+\left(C-\sum \limits_{k=1}^n x_k\right)\underbrace{p}_{=\, 0}\right\}\\
&\leqslant   \max \limits_{p \geqslant  0}\min \limits_{x_k \geqslant  0}\left\{\sum \limits_{k=1}^n f_k(x_k)+\left(C-\sum \limits_{k=1}^n x_k\right)p\right\}\\
&= \min \limits_{x_k \geqslant  0}\left\{\sum \limits_{k=1}^n f_k(x_k)+\left(C-\sum \limits_{k=1}^n x_k\right)p^{*}\right\}\\
&\leqslant   \sum \limits_{k=1}^n f_k(\bar x_k)+\left(C-\sum \limits_{k=1}^n \bar x_k\right)p^{*}\\
&\leqslant    \sum \limits_{k=1}^n f_k(\bar x_k)-Cp^{*}.
\end{split}
\end{equation*}



Отсюда получаем следующую лемму.
\begin{Lem}
\label{3_lem_pmax}
Значение $p = p^*$ (см. теорему~\ref{theo_1}), удовлетворяет неравенству:
\begin{equation*}
\label{p_location}
0 \leqslant   p^* \leqslant   p_{max},
\end{equation*}
где $p_{max}$ определяется условием 
\begin{equation}
\label{p_max}
p_{max}=\dfrac{1}{C}\sum \limits_{k=1}^n \Bigl\{f_k\left(\dfrac{2C}{n}\right)-f_k(0)\Bigr\}.
\end{equation}
\end{Lem}

\subsection{Оценка скорости сходимости}
Для обоснования сходимости метода понадобится следующая
ключевая лемма.
\begin{Lem}
\label{3_lem_eq}
\ai{Справедливо следующее неравенство
\ag{
$$f(\mathtt{x}(p))-f(\mathtt{x}^*)\leqslant  p\left(\sum \limits_{k=1}^n x_k(p) - C\right).$$}}
\end{Lem}

\textbf{Доказательство:}
Пусть $p^*$~--- решение двойственной задачи~\eqref{2_eq_dual}. Обозначим $\mathtt{x}^*=\mathtt{x}(p^*)$. 
Тогда выполняется следующее соотношение
\ag{
$$
-f(\mathtt{x}(p))-p\left(C-\sum \limits_{k=1}^n x_k(p)\right) = \varphi(p) \geqslant \varphi(p^*) = -f(\mathtt{x}^*)-p^*\underbrace{\left( C-\sum \limits_{k=1}^n x_k^*\right)}_{= \, 0}=-f(\mathtt{x}^*).
$$}
Следовательно,
\ag{$$
f(\mathtt{x}(p))-f(\mathtt{x}^*)\leqslant    p\left(\sum \limits_{k=1}^n x_k(p) - C \right).
$$}
\ag{Что и требовалось доказать.}
\qed

Поскольку $\varphi(p)$~--- гладкая функция, можно гарантировать сходимость к нулю производной: $\varphi'(p^N) \rightarrow \varphi'(p^*)=0, \; \; N \rightarrow \infty$.   


Получим оценку скорости сходимости представленного выше 
метода (см. алгоритм~\ref{alg_dih}, раздел~\ref{sect_no_f_1}). Рассмотрим следующую задачу минимизации, представляющую собой компактифицированную двойственную к исходной задаче~\eqref{2_eq_main}:
\begin{equation}
\label{3_eq_conj}
\varphi(p)\rightarrow \min \limits_{p \in Q},
\end{equation}
где $\varphi(p)$~--- выпуклая функция, и в силу леммы~\ref{3_lem_pmax},

\begin{equation}
\label{3_eq_Q}
Q=\left\{p \;| \; 0 \leqslant   p \leqslant   p_{max} \right\},
\end{equation}
где $p_{max}$ определяется из условия~\eqref{p_max}.
Липшицевость градиента двойственной задачи~\eqref{2_smooth} (см.~утверждение~1) означает, что выполнено неравенство 
\ai{
\begin{equation*}
\varphi(p^{t+1}) \leqslant   \varphi(p^{t}) + \varphi'(p^{t})(p^{t+1}-p^{t}) + \dfrac{L_{\varphi}}{2}\left(p^{t+1}-p^{t}\right)^2.
\end{equation*}
Учитывая, что $p^{t+1}=p^{t}-h\varphi'(p^{t})$, получаем
\begin{equation*}
\varphi(p^{t+1}) \leqslant   \varphi(p^{t}) - h(\varphi'(p^{t+1}))^2 + \dfrac{L_{\varphi}h^2}{2}(\varphi'(p^{t+1}))^2= \varphi(p^{t}) - h\left(1- \dfrac{L_{\varphi}h}{2}\right)(\varphi'(p^{t+1}))^2.
\end{equation*}
Отсюда, учитывая максимум по $h$ правой части достигается при $h=\dfrac{1}{L_{\varphi}}$, получаем оценку
\begin{equation*}
\label{3_eq_grad}
(\varphi'(p^N))^2\leqslant   2L_{\varphi}\Big(\varphi(p^N)-\varphi(p^{N+1})\Big)\leqslant   2L_{\varphi}\Big(\varphi(p^N)-\varphi(p^*)\Big)=2L_{\varphi}\varepsilon,
\end{equation*}
}
где $\varepsilon$~--- точность (по функции) решения двойственной задачи~\eqref{3_eq_conj}.
Используя лемму~\ref{3_lem_pmax} и лемму~\ref{3_lem_eq}, имеем
\ag{
$$
f(\mathtt{x}(p^N))-f(\mathtt{x}^*)\leqslant   p^N\left(\sum \limits_{k=1}^n x_k(p^N) - C\right) \leqslant   p^N \varphi'(p^N) \leqslant   p_{max}\left(2L_{\varphi}\varepsilon\right)^{1/2},
$$}
т.е.
\begin{equation}
\label{3_eq_6}
f(\mathtt{x}(p^N))-f(\mathtt{x}^*)\leqslant   \left(2L_{\varphi}(p_{max})^2\varepsilon\right)^{1/2}.
\end{equation}

Скорость сходимости по аргументу метода дихотомии, который используется для решения двойственной задачи, геометрическая. И, как следствие, из липшицивости функционала двойственной задачи $\varphi(p)$ (константа Липшица функционала задачи~\eqref{3_eq_conj} $M_{\varphi} \le C$ на множестве~$Q$), скорость сходимости по функции тоже геометрическая, то есть точность (по функции) решения двойственной задачи~\eqref{3_eq_conj} не будет превышать $\varepsilon$ после 
\ai{
$$
N = \log_2\left(M_{\varphi}\dfrac{\Delta p}{\varepsilon}\right)=\log_2\left(M_{\varphi}\dfrac{p_{max}}{\varepsilon}\right)
$$
итераций, где $\Delta p=p_{max}-p_{min}$~--длина начального отрезка локализации цены.}
Отсюда и из оценки~\eqref{3_eq_6} получаем, что для того, чтобы точность (по функции) решения исходной задачи~\eqref{2_eq_main} не превосходила $\varepsilon$, достаточно сделать 
$$
N = \log_2\left(\dfrac{2L_{\varphi}M_{\varphi}(p_{max})^3}{\varepsilon^2}\right)
$$
итераций метода дихотомии при решении двойственной задачи~\eqref{3_eq_conj}. Также отметим, что, так как скорость сходимости двойственной задачи~\eqref{3_eq_conj} по функции геометрическая, то скорость сходимости прямой задачи по функции тоже будет геометрической, т.е. 
$$
f(\mathtt{x}(p^N))-f(\mathtt{x}^*)\leqslant   \left(\dfrac{2L_{\varphi}M_{\varphi}(p_{max})^3}{2^N}\right)^{1/2}=\dfrac{1}{2^{N/2}}\left(2L_{\varphi}M_{\varphi}(p_{max})^3\right)^{1/2},
$$

Заметим, что решение исходной задачи с точностью $\varepsilon$ (по функции) обеспечивается следующей оценкой на ограничение
\begin{equation*}
    \left|\varphi'(p^N)\right|=\left|\sum\limits_{k=1}^n x_k(p^N)-C\right|\leqslant   \dfrac{\varepsilon}{p_{max}}.
\end{equation*}
При этом допускается, что ограничение может нарушаться, в случае, когда 
\begin{equation*}
  C -\sum\limits_{k=1}^n x_k(p^N)\leqslant   \dfrac{\varepsilon}{p_{max}},
\end{equation*}
и получаемый в результате решения задачи вектор $\mathtt{x}(p^N)=\mathtt{x}^N$ не будет принадлежать допустимому множеству.

 Оценим скорость сходимости по аргументу для прямой (исходной) задачи~\eqref{2_eq_main}. Данную оценку можно было бы получить ввиду сильной выпуклости исходной задачи~\eqref{2_eq_main}.
Однако, как было отмечено выше, получаемое значение вектора $\mathtt{x}(p^N)$ может не принадлежать допустимому множеству, так как будет нарушаться ограничение на суммарный минимальный объём производства. Поэтому, чтобы оценить сходимость по аргументу для исходной задачи, необходимо сначала спроецировать $\mathtt{x}(p^N)$ на допустимое множество. Пусть выполнено  
\begin{equation*}
  C -\sum\limits_{k=1}^n x_k(p^N) = \dfrac{\varepsilon}{p_{max}},
\end{equation*}
тогда проекцией на допустимое множество является вектор 
$$
\overline{\mathtt{x}}(p^N)=\left(x_1(p^N)+\frac{\varepsilon}{np_{max}}, \, \ldots, \, x_n(p^N)+\frac{\varepsilon}{np_{max}}\right)^{\top}.
$$
Тогда
\begin{equation*}
    \left|\left|\overline{\mathtt{x}}(p^N)-\mathtt{x}(p^N)\right|\right|_2=\dfrac{\varepsilon}{\sqrt{n}p_{max}}.
\end{equation*}
Так как $\overline{\mathtt{x}}(p^N)$  принадлежит допустимому множеству, то для него (в силу сильной выпуклости исходной задачи) выполнено 
\begin{equation*}
\label{strict_convex}
    \dfrac{\mu}{2}\left|\left|\overline{\mathtt{x}}(p^N) - \mathtt{x}^*\right|\right|_2^2 \leqslant   f(\overline{\mathtt{x}}(p^N)) - f(\mathtt{x}^*).
\end{equation*}
Учитывая липшицевость функционала исходной задачи, получаем

\begin{equation*}
\begin{split}
\left|\left|\overline{\mathtt{x}}(p^N)-\mathtt{x}(p^*)\right|\right|_2^2
&\leqslant\dfrac{2}{\mu}\left( f(\overline{\mathtt{x}}(p^N)) - f(\mathtt{x}^*)\right)\\
&\leqslant\dfrac{2}{\mu}\left( f(\mathtt{x}(p^N)) - f(\mathtt{x}^*)+M_{f}\left|\left|\overline{\mathtt{x}}(p^N)-\mathtt{x}(p^N)\right|\right|_2 \right)\\
&\leqslant\dfrac{2}{\mu}\left( \varepsilon+M_{f}\dfrac{\varepsilon}{\sqrt{n}p_{max}}\right),
\end{split}
\end{equation*}
где~$M_f$~--- константа Липшица прямой задачи~\eqref{2_eq_main}.
Тогда, по неравенству треугольника получаем 
\begin{equation*}
\begin{split}
\left|\left|\mathtt{x}(p^N)-\mathtt{x}(p^*)\right|\right|_2 
&\leqslant   \left|\left|\overline{\mathtt{x}}(p^N)-\mathtt{x}(p^*)\right|\right|_2+\left|\left|\mathtt{x}(p^N)-\overline{\mathtt{x}}(p^N)\right|\right|_2 \\
&\leqslant\left(\dfrac{2}{\mu}\left( \varepsilon +M_{f}\dfrac{\varepsilon}{\sqrt{n}p_{max}}\right)\right)^{1/2}+\dfrac{\varepsilon}{\sqrt{n}p_{max}},
\end{split}
\end{equation*}
где~$\varepsilon$~--- точность решения прямой задачи~\eqref{2_eq_main} по функции. Так как скорость сходимости по функции геометрическая, то скорость сходимости по аргументу тоже будет геометрической, т. е.
\begin{equation*}
\begin{split}
\left|\left|\mathtt{x}(p^N)-\mathtt{x}(p^*)\right|\right|_2
&\leqslant   \dfrac{1}{2^{N/4}}\left(\left(\dfrac{8L_{\varphi}M_{\varphi}(p_{max})^3}{\mu^2}\right)^{1/2}+\left(\dfrac{8L_{\varphi}M_{\varphi}M_{f}^2p_{max}}{n\mu^2}\right)^{1/2}\right)^{1/2}\\
&+\dfrac{1}{2^{N/2}}\left(\dfrac{2L_{\varphi}M_{\varphi}p_{max}}{n}\right)^{1/2}.
\end{split}
\end{equation*}

Полученный результат сформулируем в виде теоремы.
\begin{Theo}
Пусть необходимо решить задачу~\eqref{2_eq_main} в условии отсутствия у Центра информации о функциях затрат. Предполагается, что для решения этой задачи строится двойственная к ней~\eqref{3_eq_conj}, которая решается методом дихотомии. 

Тогда скорость сходимости прямой задачи по функции и по аргументу геометрическая, т.е. 
\begin{equation*}
\begin{split}
f(\mathtt{x}(p^N))-f(\mathtt{x}^*)&\leqslant\dfrac{1}{2^{N/2}}\left(\dfrac{2L_{\varphi}(p_{max})^3}{\mu}\right)^{1/2},\\
\left|\left|\mathtt{x}(p^N)-\mathtt{x}(p^*)\right|\right|_2
&\leqslant   \dfrac{1}{2^{N/4}}
\left(\left(\dfrac{8L_{\varphi}M_{\varphi}(p_{max})^3}{\mu^2}\right)^{1/2}+\left(\dfrac{8L_{\varphi}M_{\varphi}M_{f}^2p_{max}}{n\mu^2}\right)^{1/2}\right)^{1/2}\\
&+\dfrac{1}{2^{N/2}}
\left(\dfrac{2L_{\varphi}M_{\varphi}p_{max}}{n}\right)^{1/2}
.
\end{split}
\end{equation*}

В частности, чтобы точность решения исходной задачи~\eqref{2_eq_main} по функции не превышала $\varepsilon$, достаточно сделать 
$$
N = \log_2\left(\dfrac{2L_{\varphi}M_{\varphi}(p_{max})^3}{\varepsilon^2}\right)
$$
итераций метода дихотомии при решении двойственной задачи~\eqref{3_eq_conj}.
\end{Theo}
Итак, даже если функции затрат каждого предприятия в явном виде неизвестны Центру и он может только корректировать цену на рынке, суммарное количество товара, производимое предприятиями в год, сходится к необходимому количеству~$C$. Кроме того, объём производства каждого предприятия сходится к объёму производства этого предприятия, полученного в задаче~\eqref{2_eq_main}.

\section{Централизованная распределённая оптимизация}

\begin{center}
\begin{figure}[htb]
\begin{center}
\includegraphics[width=0.8\linewidth]{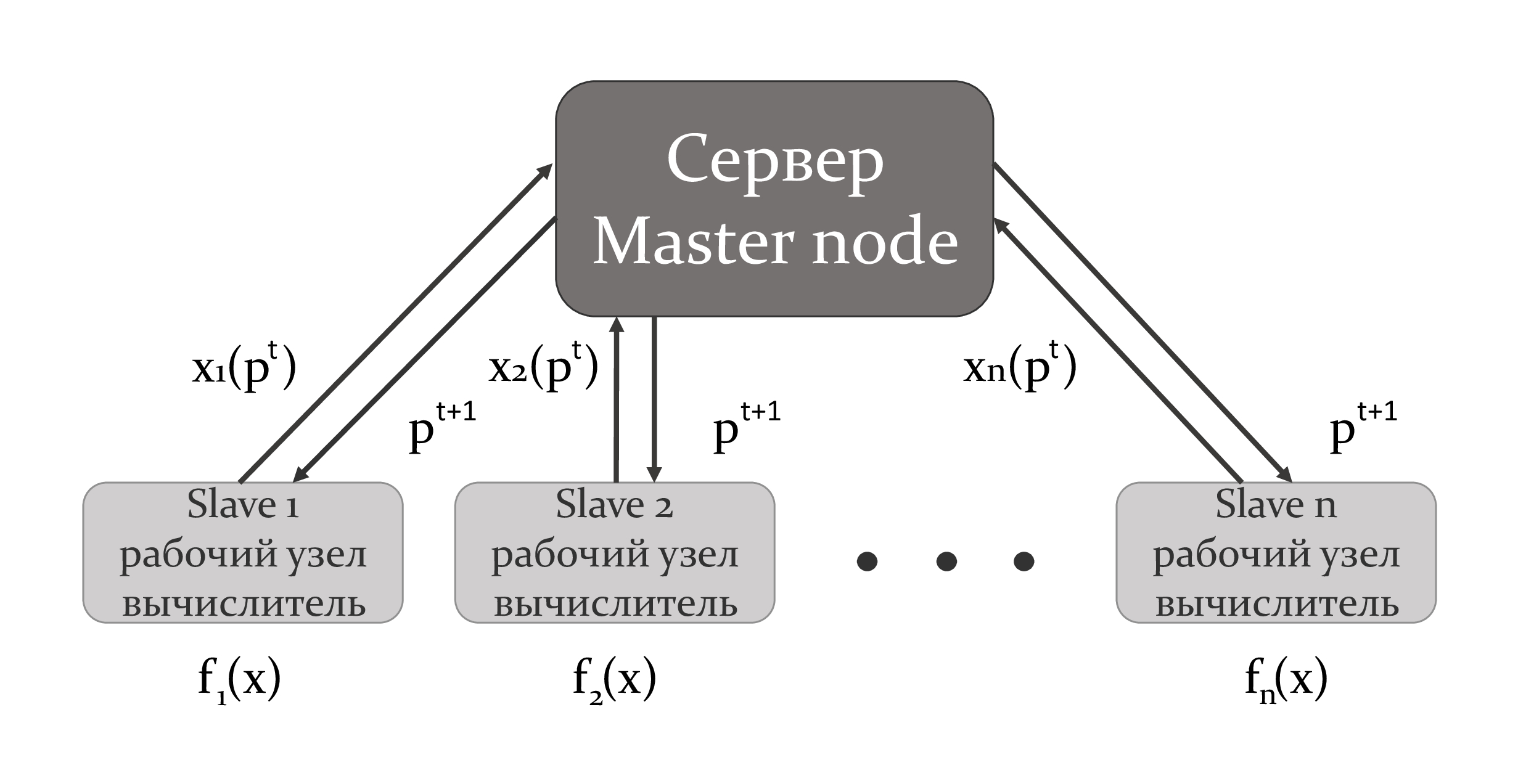}
\end{center}
\caption{Архитектура централизованной распределённой оптимизации} 
\label{arch}
\end{figure}
\end{center}

Рассмотренную в работе задачу можно проинтерпретировать и как задачу централизованной распределённой оптимизации (см., например, \cite{Nedic}). На рис.~\ref{arch} изображена архитектура централизованной распределённой оптимизации.

 Главной особенностью централизованной распределённой оптимизации является наличие одного главного сервера (Master Node), который соединяет c собой узлы (вычислители). При этом целевая функция является сепарабельной (либо её можно привести к такому виду), т.е. её можно представить в виде суммы функций, зависящих от разных компонент искомого вектора. Последнее обстоятельство позволяет организовать вычисления так, чтобы каждую компоненту вычислял отдельный узел, т.е. каждый узел обладает информацией только о своей компоненте целевой функции. Это очень удобно, так как вычисление различных частей (слагаемых) целевой функции  разными узлами происходит намного быстрее, чем если это делать на одном узле (происходит распараллеливание вычислений). Роль сервера заключается в том, чтобы посылать всем узлам значения аргументов для соответствующих компонент целевых функций и, после того как каждый узел произвёл вычисления, собрать с них результат. Далее сервер обрабатывает результат и проверяет выполнение ограничений, и на основе этого меняет значения аргументов и сообщает эту информацию вычислителям. Затем процесс повторяется снова, пока ограничения не будут выполнены и сервер не получит решение. 

 Теперь покажем связь между процессом нащупывания равновесия по Вальрасу при решении задачи~\eqref{2_eq_main} и централизованной распределённой оптимизацией. Нащупывание по Вальрасу представляет собой механизм по установлению равновесия на рынке. В основе процесса нащупывания лежит регулирование Центром (сервером) цены на абстрактный товар 
 в зависимости от разницы величин спроса и предложения на этот товар. Алгоритм централизованной распределённой оптимизации является алгоритмическим выражением нащупывания. Действительно, в нашей задаче целевая функция представляет собой суммарную функцию затрат всех предприятий, причём функция затрат каждого из предприятий зависит только от того, сколько товара данное предприятие производит. В качестве главного сервера в рассматриваемой задаче выступает Центр. А каждый вычислитель (в нашем случае предприятие) обладает только информацией о своей собственной компоненте целевой функции (непосредственно в нашей задаче это функция затрат для каждого предприятия). На текущей итерации~$t$ сервер посылает значение аргумента~$p^t$ узлам. В рассматриваемой задаче в качестве передаваемого аргумента выступает значение цены в текущем году. Далее все узлы вычисляют значения $x_k(p^t)$ (оптимальные объёмы производств при установленном Центром уровне цен). Затем сервер собирает с узлов результаты вычислений~--- $x_k(p^t)$. \ag{Затем сервер определяет знак $\sum \limits_{k=1}^n x_k(p^t) - C$ и  меняет~$p^t$ на~$p^{t+1}$. Значение $p^{t+1}$ сервер сообщает} узлам (предприятиям).


\section{Задача распределения ресурсов при отсутствии централизованного контроля цен}\label{sect_subgr}

Рассматривается задача~\eqref{2_eq_main} с неизвестными функциями затрат. На этот раз
предположим, что Центр больше не может
контролировать цены на рынке и каждое предприятие вправе устанавливать свою цену на товар. 
Пусть $k$-е предприятие установило цену~$p_k$ в текущем году и определило объём производства~$x_k$, решив
задачу максимизации прибыли. Центр сравнивает цены предприятий на товар и определяет $y_k$~--- количество тонн товара, которое будет куплено у $k$-го предприятия в текущем году (при этом Центру по-прежнему важно закупать определённое количество тонн товара в год). Далее предприятия сравнивают, на сколько и как отличается $x_k$ от $y_k$ (т.е., на сколько их <<предложение>>  отличается от <<спроса>>  Центра непосредственно на их товар), и на основе этого корректируют цену на следующий год.

Данную задачу можно переписать следующим эквивалентным образом
\begin{equation}
\label{4_eq_main}
f(\mathtt{x})=\sum \limits_{k=1}^n f_k(x_k) \rightarrow \min \limits_{\substack{\sum \limits_{k=1}^n y_k \geqslant  C, \;x_k \geqslant  y_k;\; \\ y_k \geqslant  0,\;  x_k \geqslant  0,\; k = 1, \, \ldots, \, n.}}
\end{equation}

В силу стандартной двойственности
\begin{equation*}
\begin{split}
\min \limits_{\substack{\sum \limits_{k=1}^n y_k \geqslant  C, \; x_k \geqslant  y_k,\; y_k \geqslant  0;\;\\ x_k \geqslant  0,\; k=1, \, \ldots, \, n}} f(\mathtt{x}) &= \min \limits_{\substack{\sum \limits_{k=1}^n y_k \geqslant  C\; y_k \geqslant  0;\;\\ x_k \geqslant  0,\; k=1, \, \ldots, \, n}}\Bigl\{f(\mathtt{x})+\sum \limits_{k=1}^n \max \limits_{p_k \geqslant  0}p_k(y_k-x_k)\Bigr\}\\
&=-\min\limits_{p_1, \, \ldots, \, p_n\geqslant  0} \Bigl\{\sum \limits_{k=1}^n \max \limits_{x_k \geqslant  0}(p_kx_k-f_k(x_k))-\min \limits_{\sum \limits_{k=1}^n y_k \geqslant  C; \;  y_k \geqslant  0}\sum \limits_{k=1}^{n}p_ky_k\Bigr\}\\
&=-\min\limits_{p_1, \, \ldots, \, p_n\geqslant  0}\Bigl\{\sum \limits_{k=1}^n \max \limits_{x_k \geqslant  0}(p_kx_k-f_k(x_k))-C\min \limits_{k=1, \, \ldots, \, n}p_k\Bigr\}\\
&=-\min\limits_{p_1, \, \ldots, \, p_n\geqslant  0} \Bigl\{\sum \limits_{k=1}^n \Big\{ p_kx_k(p_k)-f_k(x_k(p_k))\Big\}-C\min \limits_{k=1, \, \ldots, \, n}p_k\Bigr\},
\end{split}
\end{equation*}
где 
\begin{equation}
\label{5_eq_x_k(p)}
x_k(p_k)=\argmax \limits_{x_k \geqslant  0} \Big \{ p_k x_k-f_k(x_k) \Big \}, \quad  k =1, \, 2, \, \ldots, \, n.
\end{equation}
Тогда двойственная задача (с точностью до знака) имеет следующий вид

\begin{equation}
\label{4_eq_dual}
\varphi(p_1, \, \ldots, \, p_n)=\sum \limits_{k=1}^n \Big\{ p_kx_k(p_k)-f_k(x_k(p_k))\Big\}-C\min \limits_{k=1, \, \ldots, \, n}p_k \rightarrow \min \limits_{p_1, \, \ldots, \, p_n\geqslant  0.}
\end{equation}
Следовательно,
\begin{equation}
\label{5_eq_grad}
\nabla \varphi(p_1, \, \ldots, \, p_n)=\begin{pmatrix} 
x_1(p_1)\\ 
.\\
.\\
.\\
x_n(p_n)
\end{pmatrix} - C \begin{pmatrix} 
\lambda_1\\ 
.\\
.\\
.\\
\lambda_n
\end{pmatrix},
\end{equation}
где~$\sum \limits_{k=1}^{n} \lambda_k=1, \; \lambda_k \geqslant 0$ и если~$\lambda_k > 0$, то 
$k \, \in \, \argminset\limits_{j=1, \, \ldots, \, n} p_j$. 
Отметим, что 
здесь и в дальнейшем под~$\nabla \varphi(p_1, \, \ldots, \, p_n)$ мы будем понимать какой-то субградиент (вектор), произвольно
выбранный из выпуклого компактного множества~--- субдифференциала. При этом какой именно субградиент из конкретного субдифференциала будет выбран не является значимым, так как на оценку скорости сходимости это не повлияет. 

Под $\lambda_kC$ понимается \ag{желаемый} объём закупки Центра у~$k$-го предприятия, то есть~$\lambda_k$ определяет, какую долю от общего количества~$C$ Центр собирается закупить у~$k$-го предприятия.
В этой постановке цена никак не влияет на качество товара, т.е. все предприятия производят \ag{одинаковый товар}, только по разным ценам.

Заметим, что c экономической точки зрения субградиент~\eqref{5_eq_grad} определяет, насколько объём производства для каждого предприятия отличается от объёма желаемой закупки Центра непосредственно у этого предприятия. Учитывая это, предлагаем следующий способ определения вектора $\mathtt{\lambda}=(\lambda_1, \, \ldots, \, \lambda_n)^{\top}$.
После того, как Центр получил информацию о цене и об оптимальном объёме производства от каждого предприятия, он сравнивает полученные значения цен на товар и определяет среди них минимальное. Затем Центр распределяет желаемый суммарный объём закупки товара между предприятиями, производящими товар по (одинаковой) минимальной цене. При этом может оказаться, что суммарный объём товара по минимальной цене меньше необходимого Центру объёма. Тогда для некоторых предприятий желаемый объём закупки Центра будет превосходить произведённый ими объём товара, из чего такие предприятия смогут сделать вывод, что у Центра есть спрос непосредственно на  их товар и можно поднять цену и производить больше.     

Для решения этой задачи \eqref{4_eq_dual} предлагается использовать субградиентный спуск:\\
\begin{equation*}
\label{5_eq_subgrad_der}
\begin{pmatrix} 
p_1^{t+1}\\ 
.\\
.\\
.\\
p_n^{t+1}
\end{pmatrix}=\begin{pmatrix} 
p_1^{t}\\ 
.\\
.\\
.\\
p_n^{t}
\end{pmatrix}-h\nabla \varphi(p_1^t, \, \ldots, \, p_n^t),
\end{equation*}
где $h$~--- шаг субградиентного спуска, его значение будет
получено далее.

Заметим, что мы рассматриваем исключительно неотрицательные цены. Из-за этого наша двойственная задача становится задачей условной оптимизации, где допустимое множество задается ограничением $\mathtt{p} \geqslant  0$. Получается, что после каждой итерации метода необходимо проверять, чтобы $\mathtt{p}^{t+1}$ принадлежало допустимому множеству. Тогда метод субградиентного спуска становится методом проекции субградиента. И после вычисления нового значения $\mathtt{p}^{t+1}$ мы проектируем его на допустимое множество

\begin{equation*}
\label{5_eq_subgrad_der_1}
\mathtt{p}^{t+1}=\pi_{Q}\left(\mathtt{p}^{t}-h\nabla\varphi(\mathtt{p}^{t})\right)= \left [ \mathtt{p}^{t}-h\nabla\varphi(\mathtt{p}^{t}) \right ]_{+}.
\end{equation*} 
Непосредственно в рассматриваемой задаче проектирование на допустимое множество~--- это проверка на неотрицательность, поэтому в качестве проекции можно использовать положительную срезку\footnote{
Напомним что такое положительная срезка:
$$
\left [ x \right ]_{+} =
\begin{cases}
0, & \text{если $x\leqslant   0$;} \\
x, & \text{если $x> 0$.}
\end{cases}
$$
}.

\subsection{Стратегия поведения Центра и предприятий\\
Метод проекции субградиента}
\label{alg_2}
Опишем алгоритм поведения Центра и предприятий:

1. В текущем году $t$ (одна итерация представляет собой один год) каждое предприятие устанавливает цену на товар. Затем каждое предприятие решает оптимизационную задачу, желая максимизировать прибыль и определяет оптимальный план производства товара на текущий год, согласно заданной им цене:
$$
x_{k}(p_k^{t})=\argmax_{x_{k}\geqslant  0} \{p_k^{t}x_{k}-f_{k}(x_{k})\}, \; \;  k=1, \, \ldots, \, n,
$$
и сообщает эту информацию Центру.
Таким образом, Центр получает от предприятий вектор цен $\mathtt{p}^{t}=(p_1^{t}, \, \ldots, \, p_n^{t})^{\top}$ и вектор объёмов производства $\mathtt{x}^t(p^t)=(x_1^t(p_1^t), \, \ldots, \, x_n^t(p_n^t))^{\top}$.

2. Получив цены и объёмы производства от предприятий, Центр анализирует эту информацию и определяет долю закупок для каждого предприятия. Т.е. Центр определяет вектор $\mathtt{\lambda}^t=(\lambda^t_1, \, \ldots, \, \lambda^t_n)^{\top}$, где $\sum \limits_{k=1}^{n}\lambda^t_k=1, \; \lambda^t_k\geqslant 0$, и если $\lambda^t_k>0$, то $k \in \argminset \limits_{j=1, \, \ldots, \, n}p^t_j$. 

3. Далее каждое предприятие корректирует цену на будущий год следующим образом 
$$
\begin{pmatrix} 
p_1^{t+1}\\ 
.\\
.\\
.\\
p_n^{t+1}
\end{pmatrix}=
\left(\begin{pmatrix} 
p_1^{t}\\ 
.\\
.\\
.\\
p_n^{t}
\end{pmatrix}-h\nabla \varphi(p_1^t, \, \ldots, \,p_n^t) \right )_{+},$$
 где $$\nabla \varphi(p_1^t, \, \ldots, \,p_n^t)=\begin{pmatrix} 
x_1^t(p_1^t)\\ 
.\\
.\\
.\\
x_n^t(p_n^t)
\end{pmatrix} - C \begin{pmatrix} 
\lambda_1^t\\ 
.\\
.\\
.\\
\lambda_n^t
\end{pmatrix}.
$$
Т.е. каждое предприятие узнаёт, на сколько объём его производства отличается от желаемого объёма закупки Центра у этого предприятия в текущем году. И если Центр ничего не покупает у предприятия или покупает меньше, чем оно произвело, то предприятие понижает цену. Если же Центр готов купить больше, чем произвело предприятие, то предприятие повышает цену. В случае равенства предприятие ничего не меняет.

\subsection{Оценка скорости сходимости}

Обоснуем использование метода проекции субградиента, покажем его сходимость и оценим скорость сходимости. Стоит отметить, что используемый метод проекции субградиента является прямо-двойственным, то есть  по сходящейся (по двойственному функционалу) числовой последовательности, генерируемой методом в двойственном пространстве, можно построить последовательность в прямом пространстве, сходящуюся с той же скоростью (по функционалу в прямой задаче).

Рассмотрим прямую задачу~\eqref{4_eq_main}. Введем новую переменную $\mathtt{z}=(\mathtt{x}^{\top}, \, \mathtt{y}^{\top})^{\top} \in \mathbb{R}^{2n}$. Определим множество

\begin{equation*}
\label{4_eq_Q}
Q_{f}=\left\{\mathtt{z}=(x_1, \, \ldots, \, x_n,y_1, \, \ldots, \, y_n)^{\top} \, | \, y_k\geqslant  0, \; x_k\geqslant  0, \;k=1, \, \ldots, \,n;\;  \sum\limits_{i=1}^n y_i\geqslant  C\right\}.
\end{equation*} 
Пусть $A$~--- матрица размера \ai{$n\times 2n$}, которая имеет следующий вид:
\renewcommand{\kbldelim}{(}
\renewcommand{\kbrdelim}{)}
$$
  A = \kbordermatrix{
         & &   &   & \textbf{i}     & & &   & \textbf{i+n}   &   &\\
         & &   &   & \cdot & & &   & \cdot &   &\\
         & &   &   & \cdot & & &   & \cdot &   &\\
         & &   &   & \cdot & & &  & \cdot &   &\\
      \textbf{i}  & 0& \ldots  & 0  & -1 & 0&\ldots &0 & 1 & 0 &\ldots  & 0\\
         & &   &   &       & & &   &  & \\
         & &   &   &       & & &   &  & \\
         & &   &   &       & & &   &  &
  },
$$
т.е. в $i$-ой строке матрицы в $i$-ом столбце стоит $-1$ и в $i+n$ стоит $1$, а все остальные элементы $0$.
Тогда исходную задачу~\eqref{4_eq_main} можно представить следующим эквивалентным образом:

\begin{equation}
\label{4_eq_main_eq}
f(\mathtt{z}) \;\rightarrow \min\limits_{A\mathtt{z} \, \leqslant \,  0, \; \mathtt{z} \, \in \, Q_f}.
\end{equation}
Заметим, что данная задача не является сильно выпуклой по переменной $z$.
Двойственной задачей к задаче~\eqref{4_eq_main_eq} (с точностью до знака) будет задача:
\begin{equation}
\label{4_eq_dual_eq}
\varphi(\mathtt{p})=-\left<\mathtt{p},A\mathtt{z}(\mathtt{p})\right>-f(\mathtt{x}(\mathtt{p})) \;\rightarrow \min\limits_{p \, \geqslant  \, 0},
\end{equation}
где $\mathtt{z}(\mathtt{p})=(x_1(p_1), \, \ldots, \, x_n(p_n),y_1(p_1), \, \ldots, \, y_n(p_n))^{\top}$, $x_k(p_k),  \;k=1, \, \ldots, \,n$ определяется из условия~\eqref{5_eq_x_k(p)},  $y_k(p_k)= C\lambda_k(\mathtt{p}),\;k=1, \, \ldots, \,n$, где $\sum \limits_{k=1}^{n}\lambda_k(\mathtt{p})=1, \; \lambda_k(\mathtt{p})\geqslant 0$ и если $\lambda_k(\mathtt{p})>0$, то $k \in \argminset \limits_{j=1, \, \ldots, \, n} p_j$. В выборе $\lambda_k(\mathtt{p})$ есть произвол. Однако от того, как именно распорядиться этим произволом, не зависят получаемые далее оценки.

Пусть~$\mathtt{p}^{*}=(p^{*}_1, \, \ldots, \, p^{*}_n)^{\top}$~---  решение задачи~\eqref{4_eq_dual_eq}, а~$\varphi^*=\varphi(\mathtt{p}^{*})$~--- оптимальное значение функционала. Заметим, что в этой задаче вектор решения $\mathtt{p}^{*}$ представляет собой вектор, компоненты которого равны равновесной цене $p^*$, которую устанавливает Центр в задаче из пункта~\ref{sect_no_f}, т.е. $p^{*}_1= \, \ldots \,= p^{*}_n=p^*$. Тогда под~$\varepsilon$-приближенным решением этой задачи будем понимать такой~$\mathtt{p}^N$, что:
\begin{equation}
\label{4_eq_eps}
\varphi(\mathtt{p}^N)-\varphi^* \leqslant   \varepsilon.
\end{equation}
Согласно лемме~\ref{3_lem_pmax} значение $p^*\leqslant   p_{max}$. Тогда, исходя из определения вектора $\mathtt{p}^{*}$, получается, что он удовлетворяет следующей оценке 
\begin{equation*}
\label{4_eq_p**}
\left|\left|\mathtt{p}^{*}\right|\right|_2\leqslant   p_{max}\sqrt{n}.
\end{equation*}
Положим $\mathtt{p}^0=(0, \, \ldots, \,0)^{\top}$, где $\mathtt{p}^0$~--- вектор цен, которые предприятия устанавливают изначально.
Определим множество
$$
B^{+}_{3R}(\mathtt{p}^0) = B^{+}_{3R}(0)=\{\mathtt{p}\;:\; \mathtt{p}\geqslant  0, \; \left|\left|\mathtt{p}-\mathtt{p}^0\right|\right|_2\leqslant   3R\},
$$
где
\begin{equation}
\label{4_R}
R=\left|\left|\mathtt{p}^0-\mathtt{p}^{*}\right|\right|_2=p_{max}\sqrt{n},
\end{equation}
при этом  все получаемые $\mathtt{p}^t$ будут содержаться в $B^{+}_{2R}(\mathtt{0})$ :
\begin{equation}
\label{bound_2R}
\|\mathtt{p}^t\|_2 \leqslant 2R,
\end{equation}
поскольку (см. второй параграф книги \cite{Gasnikov})
$$\|\mathtt{p}^t\|_2 = \|\mathtt{p}^t - \mathtt{p}^0\|_2 \leqslant \|\mathtt{p}^t - \mathtt{p}^{*}\|_2 + \|\mathtt{p}^{*} - \mathtt{p}^0\|_2 \leqslant 2 \|\mathtt{p}^{*} - \mathtt{p}^0\|_2 = 2 \|\mathtt{p}^{*}\|_2 = 2R.$$
Пусть константа $M = O(C)$ определяется согласно неравенству:
$$
\left|\left|\bigtriangledown\varphi(\mathtt{p})\right|\right|_2 =\|\mathtt{x}(\mathtt{p}) - \mathtt{y}(\mathtt{p}) \|_2 =\sqrt{\sum_{k=1}^n \left(x_k(p_k) - C\lambda_k(\mathtt{p})\right)^2} \leqslant M, \; \; \; \forall \, \mathtt{p} \in B^{+}_{3R}(0).
$$
Как было описано в пункте~\ref{alg_2}, задачу условной оптимизации~\eqref{4_eq_dual_eq} предлагается решать методом проекции субградиента

\begin{equation}
\label{4_eq_pr_grad}
\mathtt{p}^{t+1}= \left [ \mathtt{p}^{t}-h\nabla\varphi(\mathtt{p}^{t}) \right ]_{+}.
\end{equation}
Из~\eqref{4_eq_pr_grad} можно получить следующее неравенство:

\begin{equation*}
\begin{split}
\left|\left|\mathtt{p}^{t+1}-\mathtt{p}^{*}\right|\right|^2_2&=\left|\left|\left(\mathtt{p}^{t}-h\nabla\varphi(\mathtt{p}^{t})-\mathtt{p}^{*}\right)_{+}\right|\right|^2_2\\
&\leqslant   \left|\left|\mathtt{p}^{t}-h\nabla\varphi(\mathtt{p}^{t})-\mathtt{p}^{*}\right|\right|^2_2\\
&=\left|\left|\mathtt{p}^{t}-\mathtt{p}^{*}\right|\right|^2_2-2h\left<\nabla\varphi(\mathtt{p}^t),\mathtt{p}^t-\mathtt{p}^{*}\right>+h^2\left|\left|\nabla\varphi(\mathtt{p}^{t})\right|\right|^2_2.
\end{split}
\end{equation*}
Перепишем его в немного другом виде (именно в таком виде мы его будем использовать в дальнейшем)

\begin{equation}
\label{4_eq_2}
\left<\nabla\varphi(\mathtt{p}^t), \, \mathtt{p}^t-\mathtt{p}^{*}\right>\leqslant   \dfrac{1}{2h}\left|\left|\mathtt{p}^{t}-\mathtt{p}^{*}\right|\right|^2_2-\dfrac{1}{2h}\left|\left|\mathtt{p}^{t+1}-\mathtt{p}^{*}\right|\right|^2_2+\dfrac{h}{2}\left|\left|\nabla\varphi(\mathtt{p}^{t})\right|\right|^2_2.
\end{equation}

В качестве критерия останова метода выберем условие

\begin{equation}
\label{4_eq_condition}
\varphi(\mathtt{p}^N)-\dfrac{1}{N}\min\limits_{\mathtt{p}\in B^{+}_{3R}(\mathtt{0})}\left\{\sum \limits_{t=0}^{N-1}\left[\varphi(\mathtt{p}^t)+\left<\nabla\varphi(\mathtt{p}^t), \, \mathtt{p}-\mathtt{p}^t\right>\right]\right\} \leqslant   \varepsilon,
\end{equation}
где $\mathtt{p}^N=\dfrac{1}{N}\sum\limits_{t=1}^{N}\mathtt{p}^t$.
При этом минимум берется по~$B^{+}_{3R}(\mathtt{0})=B^{+}_{3R}(\mathtt{p}^0)$, а не по $B^{+}_{3R}(\mathtt{p}^*)$, так как $\mathtt{p}^*$ нам не известно и, исходя из выбора $R$, $\mathtt{p}^{*}\in B^{+}_{R}(\mathtt{p}^0)$.
Заметим, что при выполнении условия~\eqref{4_eq_condition}~$\mathtt{p}^N$ удовлетворяет условию~\eqref{4_eq_eps} т.к.
\begin{equation*}
\begin{split}
\varphi(\mathtt{p}^N)-\varphi^*
&\leqslant   \varphi(\mathtt{p}^N)-\dfrac{1}{N}\sum \limits_{t=0}^{N-1}\left[\varphi(\mathtt{p}^t)+\left<\nabla\varphi(\mathtt{p}^t),\mathtt{p}^{*}-\mathtt{p}^t\right>\right]\\
&\leqslant   \varphi(\mathtt{p}^N)-\dfrac{1}{N}\min\limits_{\mathtt{p}\in B^{+}_{3R}(\mathtt{0})}\left\{\sum \limits_{t=0}^{N-1}\left[\varphi(\mathtt{p}^t)+\left<\nabla\varphi(\mathtt{p}^t),\mathtt{p}-\mathtt{p}^t\right>\right]\right\}\\
&\leqslant   \varepsilon.
\end{split}
\end{equation*}
Из неравенства~\eqref{4_eq_condition} ввиду~\eqref{4_eq_dual} получаем
\begin{center}
\begin{multline}
\varphi(\mathtt{p}^N)+\dfrac{1}{N}\sum\limits_{t=0}^{N-1}f(\mathtt{x}(\mathtt{p}^t))+\dfrac{1}{N}\sum\limits_{t=0}^{N-1}\left<\mathtt{p}^t, \, A\mathtt{z}(\mathtt{p}^t)\right>-
\dfrac{1}{N}\min\limits_{p\in B^{+}_{3R}(\mathtt{0})}\left\{\sum\limits_{t=0}^{N-1}\left<\nabla\varphi(\mathtt{p}^t), \, \mathtt{p}-\mathtt{p}^t\right>\right\}\leqslant   \varepsilon
\end{multline}
\end{center}
По выпуклости $f(\mathtt{x})$ имеем
$$
\varphi(\mathtt{p}^N)+f\Bigg( \underbrace{ \dfrac{1}{N}\sum\limits_{t=0}^{N-1}\mathtt{x}(\mathtt{p}^t)}_{\mathtt{x}^N}\Bigg) -\dfrac{1}{N}\min\limits_{\mathtt{p}\in B^{+}_{3R}(\mathtt{0})}\left\{\sum\limits_{t=0}^{N-1}-\left<A\mathtt{z}(\mathtt{p}^t),\mathtt{p}\right>\right\} =
$$
$$
\varphi(\mathtt{p}^N)+f(\mathtt{x}^N)+\dfrac{1}{N}\max\limits_{\mathtt{p}\in B^{+}_{3R}(\mathtt{0})}\left\{\sum\limits_{t=0}^{N-1}\left<A\mathtt{z}(\mathtt{p}^t), \, \mathtt{p}\right>\right\}\leqslant   \varepsilon.
$$
Подставляя в явном виде~$A\mathtt{z}(\mathtt{p}^t)=\mathtt{y}(\mathtt{p}^t)-\mathtt{x}(\mathtt{p}^t)$, получим
\begin{center}
\begin{equation*}
\varphi(\mathtt{p}^N)+f(\mathtt{x}^N)+\dfrac{1}{N}\max\limits_{\mathtt{p}\in B^{+}_{3R}(\mathtt{0})}\left\{\sum\limits_{t=0}^{N-1}\left<\mathtt{y}(\mathtt{p}^t)-\mathtt{x}(\mathtt{p}^t), \, \mathtt{p}\right>\right\}\leqslant   \varepsilon.
\end{equation*}
\end{center}
Таким образом, условие выглядит следующим образом: 
\begin{equation}
\label{4_eq_condition_2}
\varphi(\mathtt{p}^N)+f(\mathtt{x}^N)+3R\left|\left|\left(\mathtt{y}^N-\mathtt{x}^N\right)_{+}\right|\right|_{2}\leqslant   \varepsilon,
\end{equation}
где $\mathtt{y}^{N}=\dfrac{1}{N}\sum\limits_{t=0}^{N-1}\mathtt{y}(\mathtt{p}^{t})$.
Учитывая слабую двойственность~$-f(\mathtt{x}^{*})\leqslant   \varphi(\mathtt{p}^{*})$, получаем: 
\begin{equation*}
\begin{split}
f(\mathtt{x}^N)-f(\mathtt{x}^{*})&\leqslant   f(\mathtt{x}^N)+\varphi(\mathtt{p}^{*})\\
&\leqslant   \varphi(\mathtt{p}^N)+f(\mathtt{x}^N)\\
&\leqslant   \varphi(\mathtt{p}^N)+f(\mathtt{x}^N)+3R\left|\left|\left(\mathtt{y}^N-\mathtt{x}^N\right)_{+}\right|\right|_{2}\\
&\leqslant   \varepsilon 
\end{split}   
\end{equation*}
Из \eqref{4_eq_dual}, \eqref{bound_2R}, \eqref{4_eq_condition_2} получаем
 $$
 R\left|\left|\left(\mathtt{y}^N-\mathtt{x}^N\right)_{+}\right|\right|_{2}\leqslant   \overbrace{\underbrace{\dfrac{1}{N}\sum\limits_{t=0}^{N-1}\left<\mathtt{x}(\mathtt{p}^t)-\mathtt{y}(\mathtt{p}^t), \, \mathtt{p}^N\right>}_{\geqslant  -2R\left|\left|\left(\mathtt{y}^N-\mathtt{x}^N\right)_{+}\right|\right|_{2}}-f(\mathtt{x}^N)}^{\leqslant   \varphi(\mathtt{p}^N)}+f(\mathtt{x}^N)+3R\left|\left|\left(\mathtt{y}^N-\mathtt{x}^N\right)_{+}\right|\right|_{2}\leqslant   \varepsilon  
$$
Получается, что при решении задачи возникает невязка в ограничении. При этом, так как ограничение является неравенством, в невязке присутствует положительная срезка. Рассмотрим случай, когда $y_k^{N} \, \geqslant \,  x_k^{N}, \; k=1, \, \ldots, \,n$. В данном случае может нарушаться неравенство  $\sum\limits_{k=1}^{n}x_k^{N} \, \geqslant  \, C$, т.е., согласно экономической интерпретации задачи, суммарно предприятия могут произвести меньше необходимого количества $C$.
Пусть вектор $\mathtt{e}=(e_1, \, \ldots, \,e_n)^{\top}=(1, \, \ldots, \,1)^{\top}$. Тогда, используя неравенство Коши--Буняковского, получаем следующую оценку
$$\sum\limits_{k=1}^{n}y_{k}^N - \sum\limits_{k=1}^{n}x_{k}^N \leqslant \left<\mathtt{e}, \, (\mathtt{y}^N-\mathtt{x}^N)_{+}\right> \leqslant   \left|\left|\mathtt{e}\right|\right|_{2}*\left|\left|(\mathtt{y}^N-\mathtt{x}^N)_{+}\right|\right|_{2} = \sqrt{n} \|(\mathtt{y}^N - \mathtt{x}^N)_{+}\|_2.$$
Отметим, что в силу определения вектора $\mathtt{y}^{N}$ выполняется
 $\sum\limits_{k=1}^{n}y_k^N = C$. Тогда будет справедлива следующая оценка 
$$
C -\sum\limits_{k=1}^{n}x_{k}^N\leqslant   \dfrac{\varepsilon}{p_{max}}.
$$

Определим оптимальное значение для шага~$h$ субградиентного спуска:

\begin{align*}
\varphi(\mathtt{p}^N)  - \dfrac{1}{N}&\min\limits_{\mathtt{p}\in B^{+}_{3R}(\mathtt{0})}\left\{\sum \limits_{t=0}^{N-1}\left[\varphi(\mathtt{p}^t)+\left<\nabla\varphi(\mathtt{p}^t), \, \mathtt{p}-\mathtt{p}^t\right>\right]\right\}\\
&\leqslant   \dfrac{1}{N}\max\limits_{\mathtt{p}\in B^{+}_{3R}(\mathtt{0})}\left\{\sum \limits_{t=0}^{N-1}\left<\nabla\varphi(\mathtt{p}^t), \, \mathtt{p}^t-\mathtt{p}\right>\right\}
\\
&\overset{~\eqref{4_eq_2}}{\leqslant}\dfrac{1}{N}\max\limits_{\mathtt{p}\in B^{+}_{3R}(\mathtt{0})}\left\{\sum \limits_{t=0}^{N-1}\left(\dfrac{1}{2h}\left|\left|\mathtt{p}^{t}-\mathtt{p}\right|\right|^2_2-\dfrac{1}{2h}\left|\left|\mathtt{p}^{t+1}-\mathtt{p}\right|\right|^2_2+\dfrac{h}{2}\left|\left|\nabla\varphi(\mathtt{p}^{t})\right|\right|^2_2\right)\right\}
 \\
 &\leqslant   \dfrac{1}{2hN}\max\limits_{\mathtt{p}\in B^{+}_{3R}(\mathtt{0})}\left(\left|\left|\mathtt{p}^{0}-\mathtt{p}\right|\right|^2_2-\left|\left|\mathtt{p}^{N}-\mathtt{p}\right|\right|^2_2\right)+\dfrac{hM^2}{2} \\
 &\leqslant \dfrac{9R^2}{2hN}+\dfrac{hM^2}{2}\\
 &=\varepsilon.
\end{align*}
Далее необходимо минимизировать правую часть по~$h$ и тем самым получить оптимальное значение для шага субградиентного спуска:
$$
h=\argmin\limits_{x \, \geqslant \,  0}\left(\dfrac{9R^2}{2xN}+\dfrac{xM^2}{2}\right).
$$
Оптимальное значение шага субградиентного спуска
\begin{equation*}
\label{4_eq_h}
h=\dfrac{3R}{M\sqrt{N}}.
\end{equation*}
Для получения точности $\varepsilon$ необходимо сделать
\begin{equation*}
\label{4_eq_iter}
N(\varepsilon)=\dfrac{9R^{2}M^{2}}{\varepsilon^2}
\end{equation*}
итераций, где~$R$ определяется из условия~\eqref{4_R}.

Полученный результат можно сформулировать в виде следующей теоремы:

\begin{Theo}
Пусть необходимо решить задачу~\eqref{4_eq_main} в следующем смысле:
\begin{equation}
\label{4_eq_th}
f(\mathtt{x}^N)-f(\mathtt{x}^{*})\leqslant   \varepsilon,\;  C -\sum\limits_{k=1}^{n}x_{k}^N\leqslant   \dfrac{\varepsilon}{p_{max}}.
\end{equation}
Для этого рассмотрим двойственную задачу~\eqref{4_eq_dual}, которую будем решать методом проекции субградиента с точкой старта $\mathtt{p}^0=\mathtt{0}$. В качестве критерия останова выберем условия на зазор двойственности и невязку в ограничении:
$$
f(\mathtt{x}^N)+\varphi(\mathtt{p}^N)\leqslant   \varepsilon,\;  C -\sum\limits_{k=1}^{n}x_{k}^N\leqslant   \dfrac{\varepsilon}{p_{max}},
$$
из которого, в силу слабой двойственности, будет следовать выполнение нужного условия~\eqref{4_eq_th}. Тогда при выборе оптимального шага субградиентного спуска 
$$
h=\dfrac{3p_{max}\sqrt{n}}{M\sqrt{N}} = \dfrac{\varepsilon}{M^2}
$$
метод гарантированно остановится не более чем через
$$
N= \dfrac{9n(Mp_{max})^2}{\varepsilon^2}
$$
итераций.
\end{Theo}

Заметим, что аналогичную теорему можно получить со следующим способом выбора шага на  $k$-й итерации:
 $$
 h_k=\dfrac{\varepsilon}{\left|\left|\bigtriangledown\varphi(\mathtt{p}^k)\right|\right|_2^2}.
 $$

Таким образом, даже если  каждое предприятие вправе устанавливать свою цену на товар, а Центр может только выбирать, сколько закупать у каждого предприятия, суммарное количество товара, производимое предприятиями в год, сходится к необходимому количеству~$C$. Объём производства каждого предприятия сходится к объёму производства этого предприятия, полученного в задаче~\eqref{2_eq_main}.

\section{Вычислительные эксперименты}
Для проверки представленных в работе алгоритмов
они были реализованы на языке программирования~Python.

Рассмотрим сначала задачу~\eqref{2_eq_main} в постановке
п.~\ref{sect_no_f}. У каждого предприятия имеется
свой собственник, т.е. Центру функции затрат
каждого предприятия неизвестны. 
Для решения задачи был применён алгоритм~\ref{alg_dih}.
Заданная точность решения для всех задач поиска равновесной цены Центром была равна~$10^{-4}$.

\begin{figure}[htb]
\begin{subfigure}{0.33\textwidth}
    \includegraphics[width=\textwidth]{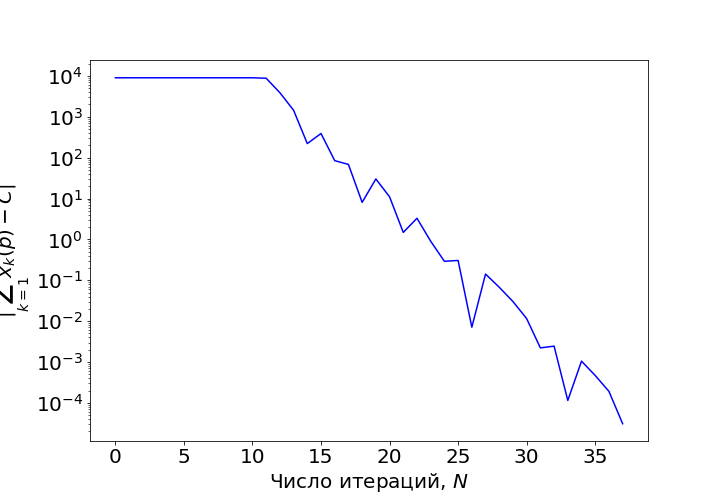}
    \caption{}
\end{subfigure}%
\begin{subfigure}{0.33\textwidth}
    \includegraphics[width=\textwidth]{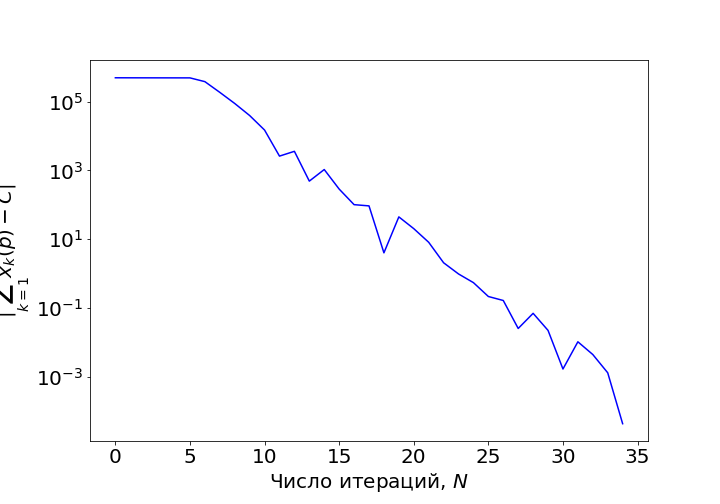}
    \caption{}
\end{subfigure}%
\begin{subfigure}{0.33\textwidth}
    \includegraphics[width=\textwidth]{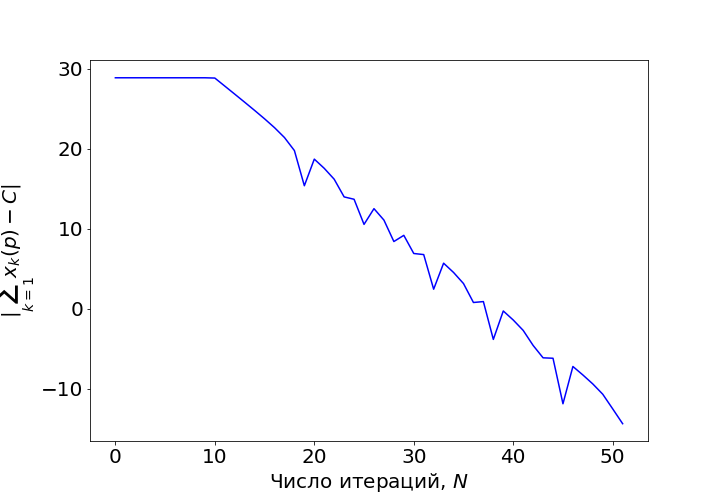}
    \caption{}
\end{subfigure}
\caption{Пример сходимости алгоритма~\ref{alg_dih} к решению для различного количества предприятий и функций затрат каждого предприятия: (a) 10 предприятий,
функции затрат~\eqref{eq_p5_prob1_fc1}; (b) 100 предприятий,
функции затрат~\eqref{eq_p5_prob1_fc21}--\eqref{eq_p5_prob1_fc22};
(c) 1000 предприятий,
функции затрат~\eqref{eq_p5_prob1_fc31}--\eqref{eq_p5_prob1_fc32}. На всех графиках масштаб по вертикальной оси логарифмический}
\label{dr_graph1}
\end{figure}

\begin{itemize}
\item Имеется~10 предприятий,
которые должны произвести~1000 условных единиц какого-то продукта. 
Цена на продукт задается Центром.
Функции затрат у всех предприятий одинаковые:
\begin{equation}
\label{eq_p5_prob1_fc1}
f_k(x_k) = \frac{1}{2} x_k^2, \quad k = 1, \, \ldots, \, 10.
\end{equation}
Решение~--- оптимальная цена~$100$~у.е.~--- было найдено за $38$~итераций. Практическая скорость сходимости показана на рис.~\ref{dr_graph1}, (a).
    \item Имеется~100 предприятий со следующими функциями затрат:
\begin{equation}
\label{eq_p5_prob1_fc21}
f_k(x_k) = \frac{1}{2} x_k^2+\frac{1}{2} x_k^4, \quad k = 1, 3 \, \ldots, \, 99,
\end{equation}
\begin{equation}
\label{eq_p5_prob1_fc22}
f_k(x_k) = 2 x_k^2, \quad k = 2, 4 \, \ldots, \, 100.
\end{equation}
Суммарно они должны произвести~$10^4$ условных единиц какого-то продукта, при этом цена на продукт задается Центром. 
Решение~--- оптимальная цена~$770.98$~у.е.~--- было найдено за $35$~итераций. Практическая скорость сходимости показана на рис.~\ref{dr_graph1}, (b).
 \item Имеется~1000 предприятий со следующими функциями затрат:
\begin{equation}
\label{eq_p5_prob1_fc31}
f_k(x_k) = x_k^2, \quad k = 1, 3 \, \ldots, \, 99,
\end{equation}
\begin{equation}
\label{eq_p5_prob1_fc32}
f_k(x_k) = 2x_k^2 + 4x_k^4, \quad k = 2, 4 \, \ldots, \, 100.
\end{equation}
Суммарно они должны произвести~$10^6$ условных единиц какого-то продукта, при этом цена на продукт задается Центром. 
Решение~--- оптимальная цена~$3987.44$~у.е.~--- было найдено за $52$~итерации.
Практическая скорость сходимости показана на рис.~\ref{dr_graph1}, (с).
\end{itemize}

\begin{figure}[htb]
\begin{subfigure}{0.5\textwidth}
    \includegraphics[width=\textwidth]{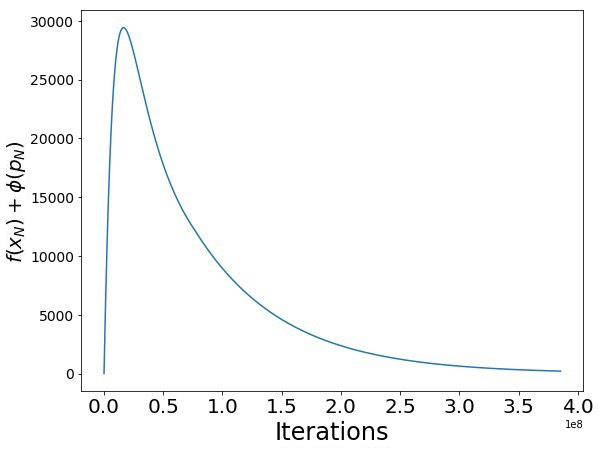}
    \caption{}
\end{subfigure}%
\begin{subfigure}{0.5\textwidth}
    \includegraphics[width=\textwidth]{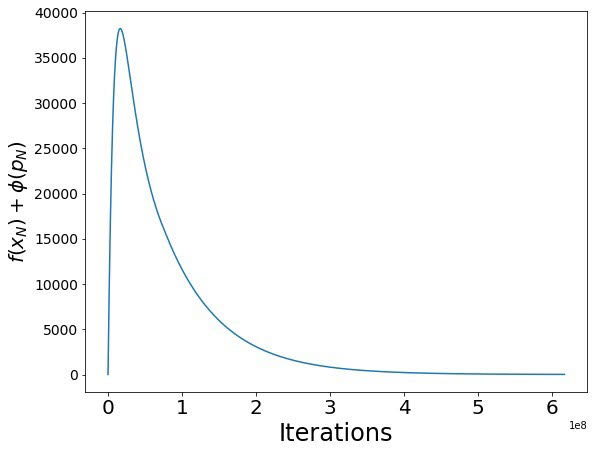}
    \caption{}
\end{subfigure}%
\caption{Пример сходимости метода проекции субградиента к решению: 
(a) заданная точность решения $\varepsilon = 10^{-1}$;
(b) заданная точность решения $\varepsilon = 10^{-2}$ }
\label{dr_graph2}
\end{figure}

Далее была численно решена задача распределения ресурсов при отсутствии централизованного контроля цен (п.~\ref{sect_subgr}).
Пусть имеются $10$~предприятий, начальный вектор цен нулевой,
точность решения $10^{-1}$ (рис.~\ref{dr_graph2}, (a)) и
$10^{-2}$ (рис.~\ref{dr_graph2}, (b)), функции затрат предприятий определяются в~\eqref{eq_p5_prob1_fc1}.
Реализация рассмотренного в п.~\ref{sect_subgr} метода решения
задачи в обоих случаях сошлась
за число итераций,
меньшее теоретической оценки.

Реализация алгоритмов на языке программирования Python выложена на Github~\cite{git}.

Работа А.В.~Гасникова, А.С.~Ивановой была поддержана грантом Президента РФ (грант МД-1320.2018.1).
Работа Е.А.~Нурминского (п.~1-3)
была поддержана грантом РФФИ 18-29-03071.
Работа Е.А.~Воронцовой в п.~1-2 была
поддержана грантом РФФИ 18-29-03071, 
в п.~4-5~--- грантом Президента РФ МД-1320.2018.1.
Авторы выражают благодарность Анжелии Недич, Ю.Е.~Нестерову \ag{и А.А.~Шананину} за ряд ценных замечаний.


\end{document}